\documentclass[11pt]{amsart}

\usepackage{amsmath,amssymb,amsthm}
\usepackage{graphicx,pxfonts}

\DeclareMathAlphabet{\mathbbold}{U}{bbold}{m}{n}

\def\k{\mathbbold{k}}

\DeclareSymbolFont{rsfscript}{OMS}{rsfs}{m}{n}
\DeclareSymbolFontAlphabet{\mathrsfs}{rsfscript}

\DeclareFontFamily{OMS}{rsfs}{\skewchar\font'177}
\DeclareFontShape{OMS}{rsfs}{m}{n}{%
      <5> rsfs5
      <6> <7> rsfs7
      <8> <9> <10> rsfs10
      <10.95> <12> <14.4> <17.28> <20.74> <24.88> rsfs10
      }{}

\def\calB{\mathrsfs{B}}

\def\calF{\mathrsfs{F}}
\def\calG{\mathrsfs{G}}
\def\calH{\mathrsfs{H}}
\def\calI{\mathrsfs{I}}
\def\calJ{\mathrsfs{J}}
\def\calK{\mathrsfs{K}}

\def\calM{\mathrsfs{M}}

\def\calO{\mathrsfs{O}}
\def\calP{\mathrsfs{P}}
\def\calQ{\mathrsfs{Q}}
\def\calR{\mathrsfs{R}}

\def\calV{\mathrsfs{V}}

\DeclareMathOperator{\Hom}{Hom}

\DeclareMathOperator{\lt}{lt}
\DeclareMathOperator{\As}{As}
\DeclareMathOperator{\Lie}{Lie}
\DeclareMathOperator{\Com}{Com}
\DeclareMathOperator{\AntiCom}{AntiCom}
\DeclareMathOperator{\Ram}{Ram}
\DeclareMathOperator{\Perm}{Perm}

\DeclareMathOperator{\Leib}{Leib}
\DeclareMathOperator{\LieGriess}{LieGriess}
\DeclareMathOperator{\PreLie}{PreLie}
\DeclareMathOperator{\End}{End}

\DeclareMathOperator{\Fin}{Fin}
\DeclareMathOperator{\Ord}{Ord}
\DeclareMathOperator{\Vect}{Vect}

\makeatletter

\theoremstyle{plain}

\newtheorem {theorem}{Theorem}

\newtheorem {corollary}{Corollary}

\newtheorem {proposition}{Proposition}

\theoremstyle{definition}

\newtheorem {definition}{Definition}
\newtheorem {remark}{Remark}
\newtheorem {example}{Example}

\begin{document}
\title{Gr\"obner bases for operads}
\author{Vladimir Dotsenko}
\address{Dublin Institute for Advanced Studies, 10 Burlington Road, Dublin 4, Ireland and School of Mathematics, Trinity College, Dublin 2, Ireland}
\email{vdots@maths.tcd.ie}
\author{Anton Khoroshkin}
\address{Departement Matematik, ETH, R\"amistrasse 101, 8092 Zurich, Switzerland and ITEP, Bolshaya Cheremushkinskaya 25, 117259, Moscow, Russia}
\email{anton.khoroshkin@math.ethz.ch}

\thanks{The first author's research was supported by grants NSh-3472.2008.2 and RFBR-CNRS-07-01-92214, and by an IRCSET research fellowship. The second 
author's research was supported by grants MK-4736.2008.1, NSh-3035.2008.2, NWO-RFBR-047017015, RFBR-07-01-00526, RFBR-CNRS-07-01-92214, and by a ETH 
research fellowship.}

\begin{abstract} We define a new monoidal category on collections (shuffle composition). Monoids in this category (shuffle operads) turn out to bring a new 
insight in the theory of symmetric operads. For this category, we develop the machinery of Gr\"obner bases for operads, and present operadic versions of 
Bergman's Diamond Lemma and Buchberger's algorithm. This machinery can be applied to study symmetric operads. In particular, we obtain an effective 
algorithmic version of Hoffbeck's PBW criterion of Koszulness for (symmetric) quadratic operads.
\end{abstract}
\maketitle

\section{Introduction}

\subsection{Description of results} Some versions of the Gr\"obner bases machinery were introduced for various algebraic structures by several authors 
(Shirshov~\cite{Shirshov} for Lie algebras, Buchberger~\cite{Buch} for commutative algebras, and Bergman~\cite{Bergman} and Bokut'~\cite{Bokut} for 
associative algebras). They proved to be extremely useful for studying various types of algebras defined by generators and relations. If one knows a 
Gr\"obner basis for relations that define a graded algebra~$A$, it can be used to compute dimensions of the graded components~$A_n$, find out whether or 
not two elements of $A$ are equal to each other etc. --- in the most efficient algorithmic way. The goal of this paper is to develop the machinery of 
Gr\"obner bases for ideals in free operads.

For nonsymmetric operads (that is, planar tree-shaped structures) some versions of Gr\"obner bases already appeared in literature, see, for example, 
\cite{Gerritzen, Hellstrom}. When one tries to define Gr\"obner bases for symmetric operads, there are two serious issues to address, both existing due to 
the symmetric groups action on components. First of all, the usual approach to Gr\"obner bases suggests that they are defined for ideals of free monoids in 
a certain monoidal category, and the definition requires a ``monomial'' basis of the free monoid together with an ordering of this basis which is 
reasonably compatible with the monoidal product. If we are working with operads, the monoidal product is given by the symmetric composition, and there is 
no known choice of an ordered basis for the free operad for which the ordering is compatible with products. Also, for the case of associative algebras a 
very important feature of Gr\"obner bases is that for every algebra they produce a monomial algebra (that is, an algebra whose defining relations are 
monomial; these algebras are usually much easier to handle than the generic ones) with the same graded dimensions and monomial basis as the original one. 
For the case of symmetric operads, it is impossible to do such a thing: together with every vanishing monomial its orbit under the symmetric group action 
should vanish, and the resulting operad will be too small (to make this reasoning precise, we suggest to the reader to look at the example of the 
operad~$\Lie$).

The crucial idea that allows to overcome the difficulties is to change the viewpoint on symmetric operads. Namely, we introduce a new monoidal structure on collections, which is different from the one used in the definition of an 
operad. We call this monoidal structure the shuffle composition of collections, and a monoid in the corresponding category~--- a shuffle operad. The word ``shuffle'' here reflects the combinatorics of compositions in this 
category; these compositions previously were somewhat implicitly used in many papers, since they provide a natural choice when computing something in an ``operad with a fixed basis''. For example, combinatorial objects very close 
to our ``shuffle surjections'' are discussed in~\cite[Prop.~11.6]{Stover}, and ``elementary shuffle compositions'' are defined in~\cite[Sect.~3.1]{Hoffbeck} (where they are called ``pointed shuffles''); see 
also~\cite[Sect.~2.2]{Merkulov}, where elementary shuffle compositions appear naturally from the geometry of face complexes for compactified configuration spaces. Shuffle operads interpolate between symmetric and nonsymmetric 
operads; they are no longer equipped with symmetric groups actions, but possible operadic compositions are not restricted to the nonsymmetric ones. The main feature of shuffle operads is that the free shuffle operad generated by 
a symmetric collection can be naturally identified with the free symmetric operad generated by the same collection. This means that we can make all necessary computations in the shuffle category, and use them to prove theorems in 
the symmetric category.

Using our methods, we define all necessary notions (orderings, reductions, $S$-polynomials) and prove an operadic analogue of Bergman's Diamond 
Lemma~\cite{Bergman}. This results in an analogue of Buchberger's algorithm~\cite{Buch} for computing Gr\"obner bases.

One of applications of Gr\"obner basis for algebras is that they give an effective version of Priddy's PBW criterion of Koszulness for quadratic 
algebras~\cite{Priddy}. Our version of Gr\"obner bases for operads provides, in a similar fashion, an effective algorithmic criterion of Hoffbeck's PBW 
criterion of Koszulness for quadratic operads~\cite{Hoffbeck}.

\subsection{Outline of the paper}

In section~\ref{OpSummary}, we discuss collections and three different monoidal structures on collections, leading to three different kinds of operads: 
nonsymmetric operads, symmetric operads, and shuffle operads. We discuss the relationship between different kinds of operads.

In section~\ref{Grobner}, we discuss combinatorics of free operads: a basis consisting of decorated trees (``tree monomials''), divisibility for monomials 
in free operads, reductions, $S$-polynomials, and Gr\"obner bases; furthermore, we present analogues of Bergman's Diamond Lemma and Buchberger's algorithm 
for computing a Gr\"obner basis. For quadratic operads, we relate our results to Hoffbeck's theory of PBW operads~\cite{Hoffbeck}.

In section~\ref{examples}, we show that for some well known operads their defining ideals have quadratic Gr\"obner bases, and present a construction which 
assigns a PBW operad to each graded commutative PBW algebra. We use this construction to prove Koszulness for some operads.

In section~\ref{further}, we discuss some possible further directions of this work, including our work in progress.

\subsection{Acknowledgements} The authors wish to thank Eric Hoffbeck for clarifying some details of~\cite{Hoffbeck} and Henrik Strohmayer for some useful 
remarks on a preliminary version of this paper. They are also grateful to Giovanni Felder, Muriel Livernet, Sergei Merkulov, Dmitri Piontkovsky, Leonid Positselski and Ivan 
Yudin for several useful discussions.

\section{Collections and monoidal structures}\label{OpSummary}
All vector spaces throughout this work are defined over an arbitrary field~$\k$ of zero characteristic. 

We denote by~$\Ord$ the category of nonempty finite ordered sets (with order-preserving bijections as morphisms), and by $\Fin$~--- the category of 
nonempty finite sets (with bijections as morphisms). Also, we denote by $\Vect$ the category of vector spaces (with linear operators as morphisms; 
unlike the first two cases, we do not require a map to be invertible).

\begin{definition}
\begin{enumerate}
 \item A \emph{(nonsymmetric) collection} is a contravariant functor from the category~$\Ord$ to the category~$\Vect$.
 \item A \emph{symmetric collection} (or a \emph{$\Sigma$-module}) is a contravariant functor from the category~$\Fin$ to the category~$\Vect$.
\end{enumerate}
For either type of collections, we can consider the category whose objects are collections of this type (and morphisms are morphisms of the corresponding 
functors).
\end{definition}

\begin{remark}
\begin{enumerate}
 \item A nonsymmetric collection is nothing but a positively graded vector space. However, the functorial definition will help us to give transparent definitions 
of monoidal structures which are otherwise (from the graded vector spaces viewpoint) totally mysterious.
 \item Let $\calP$ be a symmetric collection. Then for each finite set~$I$ the vector space~$\calP(I)$ is naturally a representation of the 
group~$\Hom_{\Fin}(I,I)$. In particular, for the ``standard'' $n$-element set $[n]=\{1,\ldots,n\}$ the vector space $\calP(n):=\calP([n])$ is a right $\Sigma_n$-module
for each~$n\ge1$. This explains the name~$\Sigma$-module.
 \item The functoriality implies that it is possible to reconstruct (in either the symmetric or the nonsymmetric case) all vector spaces $\calP(I)$ from 
the sequence of vector spaces $\calP(n)=\calP([n])$. Thus, the word ``collection'' is often used for this data.
 \item The natural forgetful functor ${}^{f}\colon\Ord\to\Fin$, $I\mapsto I^f$ leads to a forgetful functor ${}^f$ from the category of symmetric 
collections to the category of nonsymmetric ones (which literally forgets the action of the groups of symmetries):
 $$
\calP^f(I):=\calP(I^f).
 $$ 
 \item All these definitions can be given also in the case when the target category is, say, a refinement of~$\Vect$, for example, the category of graded 
vector spaces, or dg-vector spaces (chain complexes).
\end{enumerate}
\end{remark}

Now we are going to define the main ingredients used in the operad theory: monoidal structures on our categories. The first and the third one (nonsymmetric 
and symmetric compositions in the corresponding categories) are well known, the second one (the shuffle composition in the nonsymmetric category) is new. 
It provides a reasonable interpolation between the first two.

\begin{definition}
\begin{itemize}
\item Let $\calP$ and $\calQ$ be two nonsymmetric collections. Define their \emph{(nonsymmetric) composition} $\calP\circ\calQ$ by the formula
 $$
(\calP\circ\calQ)(I):=\bigoplus_{k}\calP(k)\otimes\left(\bigoplus_{f\colon I\twoheadrightarrow[k]}\calQ(f^{-1}(1))\otimes\ldots\otimes\calQ(f^{-1}(k))\right),
 $$
where the sum is taken over all non-decreasing surjections~$f$.
\item Let $\calP$ and $\calQ$ be two nonsymmetric collections. Define their \emph{shuffle composition} $\calP\circ_{sh}\calQ$ by the formula
 $$
(\calP\circ_{sh}\calQ)(I):=\bigoplus_{k}\calP(k)\otimes\left(\bigoplus_{f\colon I\twoheadrightarrow[k]}\calQ(f^{-1}(1))\otimes\ldots\otimes\calQ(f^{-1}(k))\right),
 $$
where the sum is taken over all shuffling surjections~$f$, that is surjections for which~$\min f^{-1}(i)<\min f^{-1}(j)$ whenever~$i<j$.
\item Let $\calP$ and $\calQ$ be two symmetric collections. Define their \emph{(symmetric) composition} $\calP\circ\calQ$ by the formula
 $$ (\calP\circ\calQ)(I):=\bigoplus_{k}\calP(k)\otimes_{\k S_k}\left(\bigoplus_{f\colon 
I\twoheadrightarrow[k]}\calQ(f^{-1}(1))\otimes\ldots\otimes\calQ(f^{-1}(k))\right),
 $$
where the sum is taken over all surjections~$f$.
\end{itemize}
\end{definition}

The following proposition is straightforward; we omit the proof. 

\begin{proposition} Each of the compositions defined above endows the underlying category 
with a structure of a strict monoidal category. The unit object in each case is the functor $\calI$ which vanishes on all sets of cardinality greater 
than~$1$ and is one-dimensional for any set of cardinality~$1$.
\end{proposition}

\begin{definition}
\begin{enumerate}
 \item A \emph{nonsymmetric operad} is a monoid in the category of nonsymmetric collections with the monoidal structure given by the nonsymmetric composition.
 \item A \emph{shuffle operad} is a monoid in the category of nonsymmetric collections with the monoidal structure given by the shuffle composition. 
 \item A \emph{symmetric operad} is a monoid in the category of symmetric collections with the monoidal structure given by the (symmetric) composition.
\end{enumerate}
\end{definition}

We refer the reader to~\cite{MSS} for standard background information on symmetric operads.

Note that our monoidal structures are nonlinear on the right side: the functors $\calQ\mapsto\calP\circ\calQ$ do not commute with coproducts. There is a 
convenient way to replace compositions by the so called elementary compositions which are linear in both arguments. Informally, if we interpret the 
components $\calP(n)$ of an operad $\calP$ as $n$-ary operations, the composition maps allow all possible substitutions of operations. Elementary 
compositions are those for which for all but one argument we substitute the unit element of the operad. Every composition can be described as a result of 
subsequent application of elementary compositions.

Let us describe elementary compositions more explicitly. For the sake of simplicity, we do interpret elements of operads as operations, instead of writing 
the formal categorical definitions. The proof of the proposition is omitted, as it is an immediate consequence of our definitions.

\begin{proposition}
Let $\alpha\in\calO(n)$ and $\beta\in\calO(m)$ be elements of some symmetric operad~$\calO$, and let $1\le i\le n$.
\begin{enumerate}
 \item The nonsymmetric composition $\alpha\circ_i\beta$ is the operation 
$$\alpha(x_1,\ldots,x_{i-1},\beta(x_i,x_{i+1},\ldots,x_{i+m-1}),x_{i+m},\ldots,x_{m+n-1}).$$
 \item The shuffle composition $\alpha\circ_{i,\sigma}\beta$ is the operation 
$$\alpha(x_{1},\ldots,x_{i-1},\beta(x_i,x_{\sigma(i+1)},\ldots,x_{\sigma(i+m-1)}),x_{\sigma(i+m)},\ldots,x_{\sigma(m+n-1)}).$$
Here the bijection $\sigma\colon\{i+1,\ldots, m+n-1\}\to\{i+1,\ldots, m+n-1\}$ is an $(m-1,n-i)$-shuffle, i.e.
\begin{itemize}
 \item $\sigma(i+1)<\sigma(i+2)<\ldots<\sigma(i+m-1)$,
 \item $\sigma(i+m)<\sigma(i+m+1)<\ldots<\sigma(m+n-1)$. 
\end{itemize} 
(The word ``shuffle'' reflects the way $\sigma$ permutes the elements: the relative 
order of the elements $i+1,\ldots,i+m-1$ is preserved, as well as the relative order of the elements $i+m,\ldots,m+n-1$.)
 \item The symmetric composition $\alpha\circ_{i,\sigma}\beta$ is the operation 
$$\alpha(x_{\sigma(1)},\ldots,x_{\sigma(i-1)},\beta(x_{\sigma(i)},x_{\sigma(i+1)},\ldots,x_{\sigma(i+m-1)}),x_{\sigma(i+m)},\ldots,x_{\sigma(m+n-1)}).$$ 
Here $\sigma\in S_{m+n-1}$ is an arbitrary permutation.
\end{enumerate}
\end{proposition}

For each of the three monoidal structures that we consider, one can define free monoids and (left, right, two-sided) ideals of a monoid within the general 
categorical framework (a small remark is that the categorical approach to these notions should be handled with care, since our monoidal structures do not 
commute with coproducts). We refer the reader to \cite{MSS,ValletteFree} for details. However, the categorical approach has to be translated into a working 
definition; further in this text we shall give an explicit construction for a free operad with given generators.

The main relation between our monoidal structures is described by the following 
\begin{proposition}\label{sym-shuffle}
Consider two symmetric collections $\calP$ and $\calQ$. Then we have
 $$
(\calP\circ\calQ)^f\simeq\calP^f\circ_{sh}\calQ^f.
 $$
In other words, the forgetful functor is a monoidal functor between the symmetric and the shuffle category of collections.
\end{proposition}

\begin{proof}
Consider the symmetric composition
 $$ \calP\circ\calQ(I):=\bigoplus_{k}\calP(k)\otimes_{\k S_k}\left(\bigoplus_{f\colon 
I\twoheadrightarrow[k]}\calQ(f^{-1}(1))\otimes\ldots\otimes\calQ(f^{-1}(k))\right).
 $$ 
This formula is very similar to the formula for the shuffle composition. What makes a difference is the tensor product over symmetric groups, and different 
conditions on surjections over which the direct sum is taken. It turns out that we can cover both differences simultaneously: to get rid of the symmetric group $S_k$, 
it is sufficient to fix some way to order tensor factors in the product $\calQ(f^{-1}(1))\otimes\ldots\otimes\calQ(f^{-1}(k))$, which can be, for example, 
done by introducing the condition
 $$
\min f^{-1}(1)<\min f^{-1}(2)<\ldots<\min f^{-1}(n),
 $$
which is precisely the shuffling surjection condition.
\end{proof}

There are two standard ways to define operads: via generators and relations, and via representations. Both definitions are useful, and we present them 
here. Basically, the definition via generators and relations is helpful if one wants to work inside a larger object (the free operad), and the definition 
via representations (algebras over operads) is commonly used to write down definitions and proofs, since it does not require any complicated drawings, just 
usual compositions of operations.

\begin{definition}
Fix one of the monoidal categories we are working with.
\begin{enumerate}
 \item Let $\calV$ be a collection, and let~$\calR$ be a subcollection of the free operad generated by~$\calV$. The \emph{operadic ideal} $(\calR)$ 
generated by~$\ calR$ is the minimal ideal in the free operad that contains~$\calR$. \emph{The operad with generators~$\calV$ and relations $\calR$} is the 
quotient of the free operad modulo this ideal.
 \item Let $V$ be a vector space. The \emph{operad of linear mappings} $\End_V$ is the collection $\{\End_V(n)=\Hom(V^{\otimes n},V), n\ge1\}$ of all 
multilinear mappings of $V$ into itself with the obvious composition maps.
 \item An \emph{algebra over an operad $\calO$} is a vector space~$V$ together with a morphism of the operad $\calO$ into the corresponding operad of 
linear mappings.
\end{enumerate}
\end{definition}

\begin{remark}
Consider the (either symmetric or nonsymmetric) collection $\calV$, for which $\calV(n)=0$ for $n>1$. Then the free operad generated by~$\calV$ is just the 
free associative algebra generated by~$\calV(1)$, an operadic ideal is a usual (two-sided) ideal, and an algebra over this operad (or any its quotient) is 
a (left) module over the corresponding algebra.
\end{remark}

\begin{example}
\begin{itemize}
 \item Consider the symmetric collection for which the space $\calV(2)$ is the trivial representation of~$S_2$ spanned by a binary operation~$\mu$, and all 
other spaces~$\calV(k)$ are equal to zero. The operad $\Com$ of associative commutative algebras is the maximal quotient $\calP$ of the free operad 
on~$\calV$ for which in each algebra over~$\calP$ we have $\mu(\mu(a,b),c)=\mu(a,\mu(b,c))$. It is easy to see that for all~$n$ the space $\Com(n)$ is 
one-dimensional and is spanned by the element
 $$
\mu(a_1,\mu(a_2,\mu(\ldots,\mu(a_{n-1},a_n)\ldots)).
 $$
 \item A definition of the operad $\Com$ via generators and relations: consider the free operad generated by one binary operation whose span is the trivial 
representation of $S_2$. Then the ternary component of this operad is 3-dimensional. As a representation of $S_3$, it is isomorphic to the sum of the 
trivial representation, and the 2-dimensional one. The operad $\Com$ is the quotient of the free operad modulo the ideal generated by that copy of the 
2-dimensional representation.
 \item Consider the symmetric collection for which the space $\calV(2)$ is the sign representation of~$S_2$ spanned by a binary operation~$\nu$, and all 
other spaces~$\calV(k)$ are equal to zero. The operad $\AntiCom$ is the maximal quotient $\calP$ of the free operad for which in each algebra over~$\calP$ 
we have $\nu(\nu(a,b),c)=\nu(\nu(b,c),a)$. Similarly to the case of the operad~$\Com$, one can show that the space $\AntiCom(n)$ is spanned by the 
element~$\nu(a_1,\nu(a_2,\nu(\ldots,\nu(a_{n-1},a_n)\ldots))$, but it turns out that these elements are equal to zero for~$n\ge4$ since 
\begin{multline*} 
\nu(\nu(\nu(a,b),c),d) = -\nu(\nu(a,\nu(b,c)),d) = \nu(a,\nu(\nu(b,c),d)) =\\= - \nu(a,\nu(b,\nu(c,d))) = \nu(\nu(a,b),\nu(c,d)) = -\nu(\nu(\nu(a,b),c),d). 
\end{multline*} 
Thus, the spaces $\AntiCom(n)$ are one-dimensional for $n=1,2,3$, and vanish for all other~$n$.
\end{itemize}
\end{example}

The most important consequence of Proposition~\ref{sym-shuffle} is
\begin{corollary}\label{freeop}
For a symmetric collection $\calV$, the free symmetric operad generated by $\calV$ is isomorphic, as a shuffle operad, to the free shuffle operad generated 
by $\calV^f$. Also, if $\calR$ is a symmetric subcollection of the free operad, then the ideal $(\calR)$ in the symmetric operad is isomorphic, as a 
shuffle operadic ideal, to the shuffle ideal generated by~$\calR^f$.
\end{corollary}

This leads naturally to the following idea. Assume that we wish to study a symmetric operad defined by generators and relations. As a shuffle operad, this 
operad is isomorphic to the operad defined by the same generators and relations in the shuffle category. Thus, if we have any efficient algorithms for 
working with quotients (computing bases, dimensions etc.) for shuffle operads, we can apply them, and then automatically extend the obtained results to the 
symmetric case.

\section{Gr\"obner bases}\label{Grobner}

In the case of associative algebras, to define Gr\"obner bases, one needs a monomial basis of the free algebra, and an ordering of monomials which is 
compatible with the product. In the case of operads, for a monomial basis one usually takes some class of decorated trees. Products are replaced by 
operadic compositions, which, in the symmetric case, include the action of symmetric groups on the components. To come up with a working definition of 
Gr\"obner bases, it is important to know that the monoidal structure (compositions, in our case) are reasonably compatible with the ordering of monomials. 
Instead of looking for a suitable ordered basis, we shall work with shuffle operads; according to Corollary~\ref{freeop}, we are not going to lose 
information about operadic ideals and quotients this way, in the case when our shuffle operads come from symmetric operads. Thus, we only need the 
compatibility with shuffle compositions. It turns out that for shuffle compositions there exist many different ways to define a compatible ordering. We 
shall present two of them here.

Throughout this section, the word ``operad'' means a shuffle operad. 

\subsection{Trees and a basis of the free operad}

It is well known that it is convenient to represent elements of the free operad by (decorated) trees. A (rooted) \emph{tree} is a non-empty connected directed graph $T$ of genus~$0$ for which each vertex has at least one incoming 
edge and exactly one outgoing edge. Some edges of a tree might be bounded by a vertex at one end only. Such edges are called \emph{external}. Each tree should have exactly one outgoing external edge, its \emph{output}. The 
endpoint of this edge which is a vertex of our tree is called the \emph{root} of the tree. The endpoints of ingoing external edges which are not vertices of our tree are called \emph{leaves}.

Each tree with~$n$ leaves should be labelled by~$[n]$; throughout the paper, we assume all labellings to be bijective. For each vertex $v$ of a tree, the edges going in and out of $v$ will be referred to as inputs and outputs 
at~$v$. A tree with a single vertex is called a \emph{corolla}. There is also a tree with a single input and no vertices called the \emph{degenerate} tree. Trees are originally considered as abstract graphs but to work with them 
we would need some particular representatives that we now going to describe.

For a tree with labelled leaves, its canonical planar representative is defined as follows. In general, an embedding of a (rooted) tree in the plane is determined by an ordering of inputs for each vertex. To compare two inputs of 
a vertex~$v$, we find the minimal leaves that one can reach from~$v$ via the corresponding input. The input for which the minimal leaf is smaller is considered to be less than the other one. Note that this choice of a 
representative is essentially the same one as we already made when we identified symmetric compositions with shuffle compositions.

Let us introduce an explicit realisation of the free operad generated by a collection $\calV$. The basis of this operad will be indexed by planar representative of trees with decorations of all vertices. First of all, the 
simplest possible tree is the degenerate tree; it corresponds to the unit of our operad. The second simplest type of trees is given by corollas. We shall fix a basis~$B^\calV$ of $\calV$ and decorate the vertex of each corolla 
with a basis element; for a corolla with $n$ inputs, the corresponding element should belong to the basis of~$\calV(n)$. The basis for whole free operad consists of all planar representatives of trees built from these corollas 
(explicitly, one starts with this collection of corollas, defines compositions of trees in terms of grafting, and then considers all trees obtained from corollas by iterated shuffle compositions). We shall refer to elements of 
this basis as \emph{tree monomials}.

\begin{example}\label{basis}
Let $\calO=\calF_\calV$ be the free operad for which the only nonzero component of $\calV$ is $\calV(2)$, and the basis of $\calV(2)$ is given by 
 $$
\includegraphics[scale=0.9]{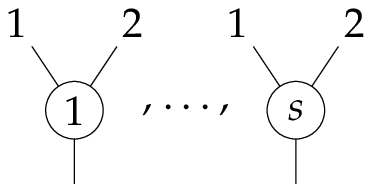}
 $$ 
Then the basis of $\calF_\calV(3)$ is given by the tree monomials
 $$
\includegraphics[scale=0.9]{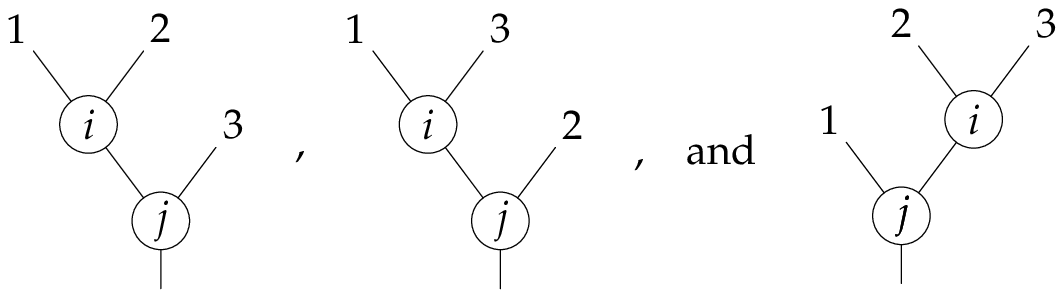}
 $$ 
with $1\le i,j\le s$.
If we assume that the $j^\text{th}$ corolla corresponds to the operation $\mu_j\colon a,b\mapsto\mu_j(a,b)$, then the above tree monomials correspond to operations
 $$
\mu_j(\mu_i(a_1,a_2),a_3), \quad \mu_j(\mu_i(a_1,a_3),a_2),\quad \text{and} \quad \mu_j(a_1,\mu_i(a_2,a_3))
 $$
respectively (as we mentioned before, this notation is much more compact than that tree notation). 
\end{example}

The free operad has two gradings that add up under compositions of operations: the arity degree, which for a tree monomial is one less than the number of leaves of the underlying tree, and the operation degree, which for a tree 
monomial is equal to the number of vertices of the underlying tree.

\begin{definition}
An element of the free operad is said to be \emph{homogeneous} if all tree monomials that occur in this element with nonzero coefficients have the same arity degree (but not necessarily the same operation degree). 
\end{definition}

\subsection{Admissible orderings}

\begin{definition}
An ordering of tree monomials of $\calF_\calV$ is said to be \emph{admissible}, if the following properties are satisfied:
\begin{itemize}
\item If $\alpha\in\calF_\calV(n)$ and $\beta\in\calF_\calV(m)$, and $n<m$, then $\alpha<\beta$.
\item If for $\alpha, \alpha'\in\calF_\calV(n)$ and $\beta,\beta'\in \calF_\calV(m)$,  we have
 $$
\alpha\le\alpha' \text{ and }\beta\le\beta',
 $$
then for all $i=1,\ldots,n$ and all $(m-1,n-i)$-shuffles $\sigma$, we have
 $$
\alpha \circ_{i,\sigma} \beta \le \alpha' \circ_{i,\sigma} \beta'.
 $$
\end{itemize}
\end{definition}

\subsubsection{Path-lexicographic ordering}

This ordering is an extension of the partial ordering defined in \cite{Hoffbeck} to a certain linear ordering.

Let $\alpha$ be a tree monomial with $n$ inputs. We associate to $\alpha$ a sequence $(a_1,a_2, \ldots,a_n)$ of $n$ words in the alphabet $B^\calV$ and a permutation $g\in S_n$ as follows. For each leaf~$i$ of the underlying tree 
$\tau$, there exists a unique path from the root to $i$. The word $a_i$ is the word composed, from left to right of the labels of the vertices of this path, starting from the root vertex. The permutation $g$ lists the labels of 
leaves of the underlying tree in the order determined by the planar structure (from left to right).

\begin{example}\label{bijection}
For the tree monomials from Example~\ref{basis}, we have 
 $$
\includegraphics[scale=0.9]{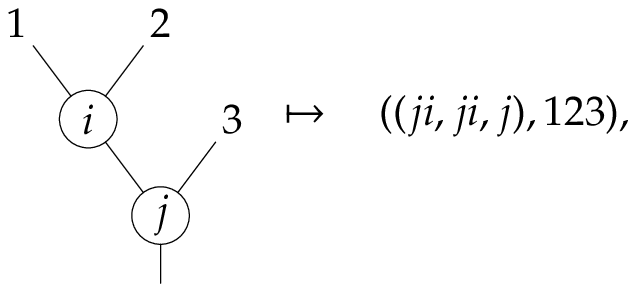} \quad\includegraphics[scale=0.9]{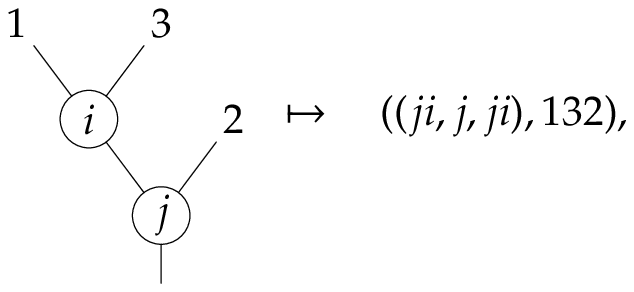}
 $$
and
 $$
\includegraphics[scale=0.9]{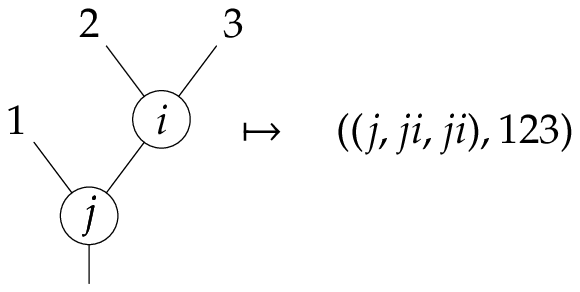}
 $$
\end{example}

The following proposition is straightforward; we leave the proof to the reader.

\begin{proposition}
This mapping from tree monomials to pairs consisting of a sequence of words and a permutation is injective. 
\end{proposition}

Now to compare two tree monomials, we just compare the corresponding sequences using the degree--lexicographic ordering: first compare the lengths of the sequences (that is, arities of our tree monomials), and if they are equal, 
compare them (word by word) using the degree--lexicographic ordering on words. If two sequences are equal to each other, we compare the permutations in reverse lexicographic order (find the first position where two permutations 
differ; the permutation for which the element at this position is smaller, is greater than the other one).

\begin{example}\label{order1}
For the tree monomials from Example~\ref{basis}, we have 
 $$
\includegraphics[scale=0.9]{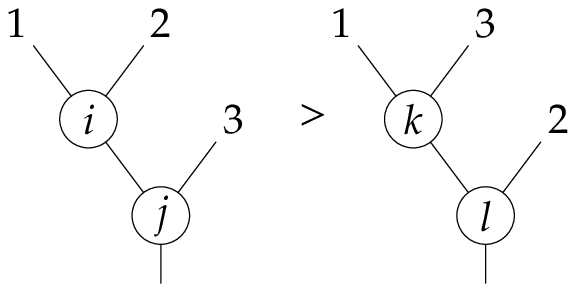}
 $$ 
if $j>l$, or $j=l$ and $i>k$, or $j=l$ and $i=k$.
\end{example}

\begin{proposition}
The ordering defined above is admissible.
\end{proposition}

\begin{proof} Let $\alpha$, $\alpha'$ be tree monomials with $n$ entries and $\beta$, $\beta'$ be tree monomials with $m$ entries. Assume that $\alpha\le\alpha'$ and $\beta\le\beta'$. 

From \cite[Prop.~3.5]{Hoffbeck}, we know that the partial order on tree monomials (ignoring the permutation data) is admissible. Thus it remains to prove our statement for the case when the underlying trees of $\alpha$ and 
$\alpha'$ are the same, and the underlying trees of $\beta$ and $\beta'$ are the same.

The composition on the level of permutations preserves all entries of the first permutation that are less than~$i$, substitutes instead of~$i$ the second permutation increased (termwise) by $i-1$, and increases all other entries 
by $m-1$. A shuffle~$\sigma$ then permutes the entries which are greater than~$i$, keeping the relative order of entries of each of the permutations unchanged.

Let $u=u_1u_2\ldots u_n$ and $u'=u'_1u'_2\ldots u'_n$ be the permutations corresponding to $\alpha$ and $\alpha'$, and $v=v_1v_2\ldots v_m$ and $v'=v'_1v'_2\ldots v'_m$ be the permutations corresponding to $\beta$ and $\beta'$. 
If $u\ne u'$, let us consider the minimal number $s$ for which $u_s<u'_s$. Let us also consider the number $r$ for which $u_r=i$, and the number $r'$ for which $u'_{r'}=i$. From the minimality of $r$, it follows that if $r<s$, 
then $r'=r$.

If $s\le r$ (so $s\le r'$ as well), then $\alpha\circ_{i,\sigma}\beta<\alpha'\circ_{i,\sigma}\beta'$, since the first $s-1$ terms of $u$ will remain equal to the first $s-1$ terms of~$u'$, and the condition $u_s<u'_s$ will also 
survive (the relative order of elements coming from the same tree is preserved by a shuffle).

If $r=r'<s$, that is, the number $i$ is on the same place in~$u$ and $u'$, then the permutations~$v$ and~$v'$ will be substituted in the same place in~$u$ and~$u'$, so $\alpha\circ_{i,\sigma}\beta<\alpha'\circ_{i,\sigma}\beta'$, 
since the relative order of elements coming from the same tree is preserved by a shuffle. It follows that either the condition $v<v'$ will work, or, if $v=v'$, the pair $u_s<u'_s$ ($m$ positions to the right compared to its 
original placement) will be the first pair where the two numbers differ.
\end{proof}

It is easy to see that this ordering admits many minor modifications: for example, to compare words we can use the lexicographic, or reverse degree-lexicographic (words of smaller degree are larger) ordering; the lexicographic 
ordering of permutations can be reversed as well.

\subsubsection{Forest-lexicographic ordering}

Consider the set $\mathsf{p}(\mathbb{N})$ of all finite subsets of $\mathbb{N}$. We define an ordering on $\mathsf{p}(\mathbb{N})$ as follows: 
 $$
I=\{i_1<\ldots<i_m\} > J=\{j_1<\ldots<j_l\}
 $$ 
if $(i_1,\ldots,i_m)>(j_1,\ldots,j_l)$ lexicographically. (Larger sets have smaller elements, for example, $[k]$ is the largest subset of cardinality~$k$.)

Let $\alpha$ be a tree monomial with inputs labelled by a finite subset~$I\subset\mathbb{N}$ (our definition will be recursive, so we prefer to not restrict ourselves to $I=[n]$.). Assume that the root of $\alpha$ is labelled by 
the generator~$A$ of the free operad, and that the root vertex has $k$ children. Thus, we have
 $$
\alpha=A(\alpha_1,\ldots,\alpha_k),
 $$
where $\alpha_j$ is a tree monomial with inputs labelled by $I_j\subset I$, $I=I_1\sqcup\ldots\sqcup I_k$, $\min I_1<\ldots<\min I_k$.
Similarly, for a tree monomial $\beta$ with inputs labelled by a finite subset~$I'\subset\mathbb{N}$ we have a similar decomposition
 $$
\beta=B(\beta_1,\ldots,\beta_l),
 $$
where $\beta_j$ is a tree monomial with inputs labelled by $I'_j\subset I'$, $I'=I'_1\sqcup\ldots\sqcup I'_l$, $\min I'_1<\ldots<\min I'_l$.

Now everything is ready for the following recursive definition.

\begin{definition}
For two tree monomials $\alpha$ and $\beta$ as above, we say that $\alpha<\beta$ if $I<I'$, or $I=I'$ and $A<B$, or $I=I'$, $A=B$, and for the smallest $k$ such that $(I_k,\alpha_k)\ne(I'_k,\beta_k)$ we have either $I_k<I'_k$ or 
$I_k=I'_k$ and $\alpha_k<\beta_k$.
\end{definition}

\begin{example}\label{order2}
For the tree monomials from Example~\ref{basis}, we have 
 $$
\includegraphics[scale=0.9]{compare.1.eps}
 $$ 
if $j>l$, or $j=l$. Thus, in general this ordering is different from the path-lexicographic ordering.
\end{example}

It is easy to prove that this ordering is compatible with shuffle compositions; we omit the proof here. Note that this definition allows many minor modifications:
for example, we can first compare the labels of the root vertices and then compare the sets of leaves.

\subsection{Divisibility in the free operad}

Take a tree monomial~$\alpha$ with the underlying tree $T$. For a subtree~$T'$ of~$T$ which contains all inputs and outputs of each its vertex, let us define a tree monomial~$\alpha'$ that corresponds to~$T'$. Its vertices are 
already decorated, so we just need to take care of the leaf labelling. For each leaf $l$ of $T'$, let us consider the smallest leaf of $T$ that can be reached from~$l$. We then number the leaves according to these ``smallest 
descendants'': the leaf with the smallest possible descendant gets the label~$1$, the second smallest~--- the label~$2$ etc.

\begin{definition}
For two tree monomials $\alpha$, $\beta$ in the free operad $\calF_\calV$, we say that \emph{$\alpha$ is divisible by $\beta$}, if there exists a subtree of the underlying tree of $\alpha$ for which the corresponding tree 
monomial $\alpha'$ is equal to~$\beta$.
\end{definition}

\begin{example} Let us give an example of divisors. We shall use the operation notation; the reconstruction of the corresponding tree monomials is left to the reader. Consider the shuffle operad generated by two binary operations 
$\alpha$ and $\beta$, and a ternary operation $\gamma$. Then the tree monomial $\mu$ corresponding to the operation
 $$
\alpha(\beta(a_1,a_3),\gamma(\beta(a_2,a_6),a_4,a_5))
 $$
has among its divisors the tree monomials corresponding to the each of operations
 $$
\alpha(\beta(a_1,a_3),a_2),\quad \alpha(a_1,\gamma(a_2,a_3,a_4)), \quad \text{ and } \gamma(\beta(a_1,a_4),a_2,a_3).
 $$
On the other hand, the tree monomial corresponding to the operation
 $$
\gamma(\beta(a_1,a_3),a_2,a_4) 
 $$
is not a divisor of $\mu$ since the ordering of its leaves does not agree with the ordering of the corresponding subtree of~$\mu$. 
\end{example}

From the fact that the free operad is generated from corollas by elementary shuffle compositions, it is easy to deduce the following

\begin{proposition}
If $\alpha$ is divisible by $\beta$, then $\alpha$ can be obtained from $\beta$ by iterations of elementary shuffle compositions with corollas. 
\end{proposition}

Assume that $\alpha$ is divisible by $\beta$. Take some sequence of compositions with corollas and elementary shuffle compositions that produces $\alpha$ from $\beta$. This sequence can be applied to any tree monomial with the 
same number of arguments as~$\beta$; we denote that operation on tree monomials by $m_{\alpha,\beta}$. It is easy to see that this operation is actually well defined (that is, depends only on~$\alpha$ and its subtree 
corresponding to the divisor~$\beta$, but not on the sequence of compositions that create $\alpha$ from~$\beta$). Note that by the construction $m_{\alpha,\beta}(\beta)=\alpha$, and from the fundamental property of operadic 
orderings it is clear that if $\gamma<\beta$, then $m_{\alpha,\beta}(\gamma)<\alpha$.

\subsection{Reductions and $S$-polynomials}

All results of this section and the further ones are valid for every admissible ordering of tree monomials. 

\begin{definition}
For an element~$f$ of the free operad, the tree monomial~$\alpha$ is said to be its \emph{leading term}, if it is the largest of the terms which occur in the expansion of $f$ with a nonzero coefficient (notation: 
$\lt(f)=\alpha$). This nonzero coefficient (the leading coefficient of $f$) is denoted by $c_f$.
\end{definition}

\begin{definition}
Assume that $f$ and $g$ are two homogeneous elements of $\calF_\calV$ for which the leading term of $f$ is divisible by the leading term of $g$. The element
 $$
r_g(f):=f-\frac{c_f}{c_g}m_{\lt(f),\lt(g)}(g),
 $$
is called the \emph{reduction of $f$ modulo~$g$}. Note that by the construction we have $\lt(r_g(f))<\lt(f)$.
\end{definition}

\begin{definition}
A tree monomial $\gamma$ is called a \emph{common multiple} of two tree monomials $\alpha$ and $\beta$, if it is divisible by both $\alpha$ and $\beta$. Tree monomials $\alpha$ and $\beta$ are said to have a \emph{small common 
multiple}, if they have a common multiple for which the number of vertices of the underlying tree is less than the total number of vertices for $\alpha$ and $\beta$.
\end{definition}

\begin{definition}
Assume that $f$ and $g$ are two homogeneous elements of $\calF_\calV$ whose leading terms have a small common multiple $\gamma$. We have
 $$
m_{\gamma,\lt(f)}(\lt(f))=\gamma=m_{\gamma,\lt(g)}(\lt(g)).
 $$
The element
 $$
s_\gamma(f,g):=m_{\gamma,\lt(f)}(f)-\frac{c_f}{c_g} m_{\gamma,\lt(g)}(g),
 $$
is called the $S$-polynomial of $f$ and $g$ (corresponding to the common multiple $\gamma$; note that there can be several different small common multiples). 
\end{definition}

\begin{remark}
$S$-polynomials, as defined here, include the reductions as a particular case. It turns out to be convenient, but we shall need reductions on their own to deal with Gr\"obner bases. 
\end{remark}

\subsection{Gr\"obner bases}

In this section and further on, we assume that $\calM$ is an operadic ideal of~$\calF_\calV$, and $\calG$ is a system of homogeneous generators of~$\calM$. 

\begin{definition}
$\calG$ is called a \emph{Gr\"obner basis} of~$\calM$, if for every $f\in\calM$ the leading term of $f$ is divisible by the leading term of some element of~$\calG$.
\end{definition}

\begin{example}
In this example, we use path-lexicographic ordering.
\begin{enumerate}
 \item For the operad $\Com$, the space of generators is one-dimensional, so no additional ordering is required. We have 
 $$
\mu(\mu(a_1,a_2),a_3) > \mu(\mu(a_1,a_3),a_2) > \mu(a_1,\mu(a_2,a_3)),
 $$
and the elements 
 $$
g_1=\mu(\mu((a_1,a_2),a_3)-\mu(a_1,\mu(a_2,a_3))\text{ and }g_2=\mu(\mu(a_1,a_3),a_2) -\mu(a_1,\mu(a_2,a_3))
 $$ 
form a Gr\"obner basis of the ideal that defines $\Com$. Indeed, we know that the elements $b_n=\mu(a_1,\mu(a_2,\mu(\ldots,\mu(a_{n-1},a_n)\ldots))$ form a basis of $\Com$. These elements in each arity are the smallest monomials 
of the corresponding arity, so the leading terms of the elements of the defining ideal are all the remaining tree monomials. As one can easily see, every tree monomial which is different from one of~$b_n$ can be reduced modulo 
$\{g_1,g_2\}$. It follows that $g_1$ and $g_2$ form a Gr\"obner basis.
 \item In the case of the operad $\AntiCom$, we know that besides the generators 
 $$
g_1=\nu(\nu(a_1,a_2),a_3)+\nu(a_1,\nu(a_2,a_3)) \text{ and } g_2=\nu(\nu(a_1,a_3),a_2)-\nu(a_1,\nu(a_2,a_3)),
 $$
the element 
 $$
g_3=\nu(a_1,\nu(a_2,\nu(a_3,a_4)))
 $$
is also equal to zero in the quotient. Thus, the only difference from the operad~$\Com$ is that starting from arity~$4$, all tree monomials belong to the basis of the ideal. All of them except for the smallest one can be 
reduced modulo~$\{g_1,g_2\}$, and the smallest one can be reduced modulo~$g_3$. It follows that~$\{g_1,g_2,g_3\}$ is a Gr\"obner basis.
\end{enumerate}
\end{example}

\begin{definition}
The element $f\in \calF_\calV$ is said to have \emph{the residue $\overline{f}$ modulo $\calG$}, if $f-\overline{f}\in\calM$, and $\overline{f}$ is a linear combination of tree monomials none of which have nontrivial reductions 
modulo~$\calG$.  Notation: $f\equiv\overline{f}\pmod{\calG}$.
\end{definition}

Generally the residue modulo~$\calG$ is not unique. Nevertheless it is unique if $\calG$ is a Gr\"obner basis.

\begin{proposition}
If $\calG$ is a Gr\"obner basis of~$\calM$, the residue of each element $f\in\calF_\calV$ modulo~$\calG$ is well defined.
\end{proposition}

\begin{proof}
If $f_1$ and $f_2$ are two different residues, we have $f-f_1\in\calM$ and $f-f_2\in\calM$, so $f_1-f_2\in\calM$. The element $f_1-f_2$ is a non-zero combination of monomials that cannot be reduced modulo~$\calG$, which is a 
contradiction since the leading term of this element should be divisible by the leading term of some element of $\calG$.
\end{proof}

\begin{definition}
If $\calG$ is a Gr\"obner basis of~$\calM$, the residue of an element $f\in\calF_\calV$ is called the \emph{normal form} of~$f$. 
\end{definition}

\begin{corollary}
Let $\calG$ be a Gr\"obner basis for $\calM$. Then tree monomials that cannot be reduced modulo~$\calG$ form a basis of the quotient~$\calF_\calV/\calM$.
\end{corollary}

\subsection{Diamond Lemma}

\begin{definition}
For an element $f\in\calM$, an expansion
 $$
f=f_1+\ldots+f_m,
 $$
where the elements $f_i$ are obtained from some elements $g_i\in\calG$ by shuffle compositions, 
is said to be an \emph{admissible $\calG$-representation}, if 
 $$
\lt(f)=\max\lt(f_k).
 $$
\end{definition}

\begin{theorem}\label{diamond}
The following properties are equivalent:
\begin{itemize}
 \item[(i)] $\calG$ is a Gr\"obner basis for $\calM$;
 \item[(ii)] for all $f\in\calM$, we have $f\equiv 0 \pmod{\calG}$;
 \item[(iii)] for all pairs of elements from~$\calG$, all their $S$-polynomials (if defined) are congruent to zero modulo~$\calG$.
 \item[(iv)] every $f\in\calM$ has an admissible $\calG$-representation.
\end{itemize}
\end{theorem}

\begin{proof}

(i) implies (ii): indeed, let us use induction on the leading monomial of $f$. If $f\ne0$, then its leading term is divisible by the leading term of some $g\in\calG$. The reduction of $f$ modulo~$g$ belongs to $\calM$ and its 
leading monomial is less than the leading monomial of~$f$, so the induction hypothesis applies.

(ii) implies (iii): indeed, all reductions and $S$-polynomials belong to $\calM$, so they are congruent to zero modulo~$\calG$, because all elements of~$\calM$ are.

(iii) implies (iv): Let us prove that for any representation of $f$ as a combination of elements obtained from elements of $\calG$ by shuffle compositions, if we have $\lt(f)<\max(\lt(f_i))$, then we can find another 
representation $f=f_1'+\ldots+f_{m'}'$ for which $\max(\lt(g_k'))<\max(\lt(g_k))$. Without loss of generality, we can assume that in our representation of~$f$ each summand $f_i$ is a scalar multiple of an element obtained from a 
certain element $g_i\in\calG$ by iterated shuffle compositions with tree monomials (that is, is a ``monomial multiple'' of $g_i$).

In order for the leading term of the result to be less than the maximal leading term, some leading terms have to cancel. We shall prove our statement by induction on~$N$, the number of indices $j$ for which 
$\lt(f_j)=\max(\lt(f_k))$. Since the maximal leading term have to cancel, we have $N\ge2$. Without the loss of generality, we assume that $\max(\lt(f_k))=\lt(f_1)=\lt(f_2)$. This, in turn, means that $\lt(f_1)$ is divisible by 
both $\lt(g_1)$ and $\lt(g_2)$. Let us consider two different cases: (1) the underlying trees of these leading terms have common edges or (2) these underlying trees are disjoint.

In the first case, the leading terms of $g_1$ and $g_2$ have a small common multiple~$g$, and the leading terms of $f_1$ and $f_2$ are divisible by~$g$. Thus, 
 $$
f_1=c_1 m_{\lt(f_1),\lt(g_1)}(g_1)=c_1 m_{\lt(f_1),g}(m_{g,\lt(g_1)}(g_1))
 $$
and
 $$
f_2=c_2 m_{\lt(f_2),\lt(g_2)}(g_2)=c_2 m_{\lt(f_2),g}(m_{g,\lt(g_2)}(g_2))
 $$
for some constants $c_1$ and $c_2$.
Note that by the definition of $S$-polynomials we have
 $$
m_{g,\lt(g_1)}(g_1)=s_g(g_1,g_2)+\frac{c_{g_1}}{c_{g_2}} m_{g,\lt(g_2)}(g_2),
 $$
so (here we also take into account that $\lt(f_1)=\lt(f_2)$)
\begin{multline*}
f_1+f_2=c_1 m_{\lt(f_1),g}(m_{g,\lt(g_1)}(g_1))+f_2=\\
=c_1 m_{\lt(f_2),g}(s_g(g_1,g_2)+\frac{c_{g_1}}{c_{g_2}} m_{g,\lt(g_2)}(g_2))+f_2=\\
=\frac{c_1}{c_2}\frac{c_{g_1}}{c_{g_2}} \left(c_2 m_{\lt(f_2),g}(m_{g,\lt(g_2)}(g_2))\right)+c_1 m_{\lt(f_1),g}(s_g(g_1,g_2))+f_2=\\
=\left(1+\frac{c_1}{c_2}\frac{c_{g_1}}{c_{g_2}}\right)f_2+c_1 m_{\lt(f_1),g}(s_g(g_1,g_2)).
\end{multline*}
The leading term of the second summand is less than the maximal leading term in our representation. Thus, the number of summands with the leading term equal to $\max(\lt(f_k))$ has decreased, and we can use the induction hypothesis.

In the second case, the leading terms of $g_1$ and $g_2$ are disjoint in 
 $$
\max(\lt(f_k))=\lt(f_1)=\lt(f_2).
 $$ 
This means that (up to a scalar multiple) $f_1$ can be obtained from $g_1$ by shuffle compositions in such a way that at some point we compute the shuffle composition with the leading term of $g_2$. Let us replace this occurence 
of that leading term by $(\lt(g_2)-g_2)+g_2$. As a result, we represent $f_1$ as a sum of an element~$f_1'$ which is obtained from $\lt(f_1)$ by replacing the occurences of $\lt(g_1)$ and $\lt(g_2)$ by $g_1$ and $g_2$ 
respectively, and a certain element from $\calM$ whose leading term is less than~$\max(\lt(f_i))$. The same works for~$f_2$, and the elements $f_1'$ and $f_2'$ are proportional. Thus we can join them together, decreasing the 
number of summands with the leading term equal to $\max(\lt(f_k))$, and the induction hypothesis applies.

(iv) implies (i): indeed, for an admissible representation of $f$, let us take $f_l$ for which $\lt(f)=\lt(f_l)$. Since $f_l$ is obtained from $g_l$ by a sequence of shuffle compositions, the leading term of~$f_l$ is divisible by 
the leading term of~$g_l$.
\end{proof}

\subsection{Buchberger's Algorithm}

From our results, we immediately obtain an analogue of the Buchberger's algorithm for operads. The input of the algorithm is a set~$\calK$ of generators for the ideal~$\calM$. The algorithm works as follows

\begin{itemize}
 \item[(i)] Compute all reductions of elements from~$\calK$ modulo each other, and all pairwise $S$-polynomials. Reduce all these elements modulo~$\calK$ until they cannot be reduced further. Extend~$\calK$ by joining these 
reductions to it.
 \item[(ii)] Repeat the step (i) until there are no nonzero elements joined.
\end{itemize}

\begin{remark}
For the case of symmetric operads, we should add one more step in the very beginning: 
\begin{quote}
 Extend $\calK$ by all the results of the symmetric group actions on it, so that the subspace spanned by~$\calK$ is stable under the symmetric group action. 
\end{quote}
Then we can guarantee that the shuffle ideal generated by~$\calK$ is isomorphic, as a shuffle ideal, to the symmetric ideal generated by~$\calK$, which is sufficient for our purposes.
\end{remark}

\begin{theorem}
The output of our algorithm results is a (possibly infinite) Gr\"obner basis for the ideal~$\calJ$.
\end{theorem}

\begin{proof}
By the construction, all reductions and $S$-polynomials for the result of our algorithm can be reduced to zero, hence by our criterion, this result is indeed a Gr\"obner basis.
\end{proof}

\begin{remark}
An operad generated by unary operations is just an associative algebra. In this case, Theorem~\ref{diamond} is precisely Bergman's Diamond Lemma~\cite{Bergman} and our previous algorithm is precisely Buchberger's algorithm~\cite{Buch}.
\end{remark}

\begin{definition}
A Gr\"obner basis~$\calG$ of~$\calM$ is said to be \emph{reduced}, if there are no reductions of its elements modulo each other, or, in other words, if there are no pairs of elements of~$\calG$ for which the leading term of one 
of them is divisible by the leading term of the other one.
\end{definition}

Reduced Gr\"obner bases are useful in practice, since in general we would prefer a Gr\"obner basis to be as small as possible (the leading terms of a Gr\"obner basis give the list of tree monomials whose multiples are not allowed 
in the basis, so the problem of determining basis elements is easier if this list is shorter).

\begin{remark}
The above algorithm can be easily modified so that its output is a reduced Gr\"obner basis. Indeed, when computing reductions, we can remove from~$\calK$ the element that is being reduced, replacing it by its reduction.
\end{remark}

In all the examples of computations throughout this paper, we use the path-lexicographic ordering. When we write down linear combinations of tree monomials, we underline leading terms to make our computations easier to follow.

\begin{example}
Let us show how our algorithm applies to the two simplest examples for which we already know Gr\"obner bases.
\begin{enumerate}
 \item For the operad $\Com$, the elements 
 $$
\underline{\mu(\mu((a_1,a_2),a_3)}-\mu(a_1,\mu(a_2,a_3))\text{ and }\underline{\mu(\mu(a_1,a_3),a_2)} -\mu(a_1,\mu(a_2,a_3))
 $$ form a Gr\"obner basis of the ideal that defines $\Com$. Indeed, there are no nontrivial reductions, the action of the symmetric group preserves the subspace spanned by these elements, and all the $S$-polynomials can be reduced to zero; for example, the element $\mu(\mu(\mu(a_1,a_3),a_2),a_4)$ which is divisible by both leading terms gives rise to the $S$-polynomial
 $$
S=\mu(\mu(a_1,a_3),\mu(a_2,a_4))-\underline{\mu(\mu((a_1,\mu(a_2,a_3)),a_4)},
 $$
for which the sequence of reductions is as follows:
\begin{multline*}
S\mapsto\underline{\mu(\mu(a_1,a_3),\mu(a_2,a_4))}-\mu(a_1,\mu(\mu(a_2,a_3),a_4))\mapsto \\
\mapsto -\underline{\mu(a_1,\mu(\mu(a_2,a_3),a_4))} +\mu(a_1,\mu(\mu(a_2,a_4),a_3))\mapsto \\
\mapsto \underline{\mu(a_1,\mu(\mu(a_2,a_4),a_3))} - \mu(a_1,\mu(a_2,\mu(a_3,a_4)))\mapsto 0.
\end{multline*}
This gives yet another proof of the fact that the elements $$\mu(a_1,\mu(a_2,\mu(\ldots,\mu(a_{n-1},a_n)\ldots))$$ form a basis of $\Com$ (they are the only monomials that cannot be reduced modulo the leading terms).
 \item In the case of the operad $\AntiCom$, we start with the generators 
 $$
\underline{\nu(\nu(a_1,a_2),a_3)}+\nu(a_1,\nu(a_2,a_3)) \text{ and } \underline{\nu(\nu(a_1,a_3),a_2)}-\nu(a_1,\nu(a_2,a_3))
 $$
of degree~$2$ in operations. The common multiple $\nu(\nu(\nu(a_1,a_4),a_2),a_3)$ of the leading terms gives us an $S$-polynomial 
 $$
S=\nu(\nu(a_1,a_4),\nu(a_2,a_3))+\underline{\nu(\nu(a_1,\nu(a_2,a_4)),a_3)}, 
 $$
for which the sequence of reductions is as follows:
\begin{multline*}
S\mapsto \underline{\nu(\nu(a_1,a_4),\nu(a_2,a_3))}-\nu(a_1,\nu(\nu(a_2,a_4),a_3))\mapsto\\
\mapsto -\underline{\nu(a_1,\nu(\nu(a_2,a_3),a_4))}+\nu(a_1,\nu(\nu(a_2,a_4),a_3))\mapsto \\ 
\mapsto \nu(a_1,\nu(a_2,\nu(a_3,a_4)))+\underline{\nu(a_1,\nu(\nu(a_2,a_4),a_3))}\mapsto  2\nu(a_1,\nu(a_2,\nu(a_3,a_4))),
\end{multline*}
which cannot be reduced further. Thus we have to add the element $$\nu(a_1,\nu(a_2,\nu(a_3,a_4)))$$ to our Gr\"obner basis. Instead of computing other $S$-polynomials, let us notice that two leading terms of the quadratic relations already give the quotient of dimension~$1$ in each arity, and now our new relation kills all the components of arity at least $4$, so any further elements (of higher arities) that we might get cannot really put any new restrictions, so there is no need to join anything else to get a Gr\"obner basis.
\end{enumerate}
\end{example}

\subsection{Another Gr\"obner basis criterion}\label{GrobnerCrit}

In this section, we prove yet another criterion for a set of elements in an operadic ideal to be a Gr\"obner basis; unlike the one we proved before, this criterion is hard to check, but it is useful for theoretical applications of Gr\"obner bases. For quadratic operads (and $k=2$ in the third condition below), this definition coincides with the definition of PBW bases given by Hoffbeck~\cite{Hoffbeck} (up to reversing the ordering: the definition of PBW operads in \cite{Hoffbeck} represents non-basis compositions as combinations of larger terms). Thus, a PBW operad in the sense of Hoffbeck is just an operad with a quadratic Gr\"obner basis. 

\begin{definition}
Let $\calP$ be an operad, $\calP\simeq\calF_\calV/\calM$. A set of tree monomials $B^\calP\supset\calV$ in the free operad $\calF_\calV$ is said to be a \emph{$k$-triangular basis} of~$\calP$ if 
\begin{enumerate}
\item The image of $B^\calP$ under the canonical projection $$\calF_\calV\twoheadrightarrow \calF_\calV/\calM\simeq\calP$$ is a basis of $\calP$.
\item For $\alpha,\beta\in B^\calP$, and $\sigma$ a shuffle, either $\alpha \circ_{i,\sigma} \beta$ is in $B^\calP$, or the elements of the basis $\gamma \in B^\calP$ which appear in the expansion 
 $$
\alpha \circ_{i,\sigma} \beta=\sum_\gamma c_\gamma \gamma+\calM
 $$ 
satisfy $\gamma < \alpha \circ_{i,\sigma} \beta$ in $\calF_\calV$.
\item A tree monomial $\alpha$ belongs to $B^\calP$ if and only if for every its subtree with at most $k$ vertices the corresponding restricted tree monomial belongs to $B^\calP$.	
\end{enumerate}
\end{definition}

\begin{theorem}
Let $\calP\simeq\calF_\calV/\calM$ be an operad, $\calG\subset\calM$ be a system of homogeneous generators. Then if $\calG$ is a Gr\"obner basis for $\calM$, then the set of tree monomials which are not divisible by leading terms of $\calG$ is a $k$-triangular basis of $\calP$, where $k$ is the maximal number of vertices in leading terms of elements of~$\calG$. Conversely, for any $k$-triangular basis~$\calB$, there exists a Gr\"obner basis whose elements are combinations of tree monomials with at most~$k$ vertices which produces $\calB$ in the way described above.
\end{theorem}

\begin{proof}
Assume that $\calG$ is a Gr\"obner basis. Then the first and the third $k$-triangular basis conditions are satisfied automatically, and the second condition is satisfied as well, since the composition of two basis elements is brought to its normal form via reductions that lower the leading term.

Conversely, assume that $\calP$ has a triangular basis. For every two tree monomials $\alpha$ and $\beta$ with at most $k$ vertices (in total) and every shuffle~$\sigma$ such that~$\alpha\circ_{i,\sigma}\beta$ does not belong to the triangular basis, consider the expression 
 $$
\alpha\circ_{i,\sigma}\beta=\sum_{\gamma<\alpha\circ_{i,\sigma}\beta} c_\gamma \gamma +\calM,
 $$ 
which exists by definition of the triangular basis. Let
 $$
h_{\alpha,\beta,\sigma}:=\alpha\circ_{i,\sigma}\beta-\sum_{\gamma<\alpha\circ_{i,\sigma}\beta} c_\gamma \gamma.
 $$
Denote by $\calH$ the system of elements $\{h_{\alpha,\beta,\sigma}\}$.  It is easy to see that every tree monomial $\delta$ that does not belong to $B^\calP$ can be reduced modulo $\calH$ to a linear combination of elements of $B^\calP$ which are strictly less than~$\delta$. Indeed, because of the third triangular basis condition we can find a subtree with at most~$k$ vertices which does belong to $B^\calP$ but has its all proper subtrees in $B^\calP$, and apply the second triangular basis condition to this subtree, replacing it by a combination of smaller elements, which allows us to continue by induction. Hence every element $f\in\calF_\calV$ can be reduced to a linear combination of elements of $B^\calP$ which are strictly less than the leading term of~$f$. For an element of the ideal~$\calM$, this means that it can be reduced to zero (since the projection from~$\calF_\calV$ to~$\calP$ is injective on~$B^\calP$), which is the definition of the Gr\"obner basis. 
\end{proof}

From \cite[Prop.~3.10]{Hoffbeck}, we get the following

\begin{corollary}
An operad with a quadratic Gr\"obner basis is Koszul. 
\end{corollary}

This result seems to cover most of widely used examples of Koszul operads. We shall discuss some of these examples below.

\section{Examples}\label{examples}

\subsection{Computation of Gr\"obner bases}

In this section, we compute Gr\"obner bases for some well known operads, so that the reader can see how our machinery actually works. We concentrate on examples of operads whose Gr\"obner bases are quadratic. A good reason to do so is because the existence of a quadratic Gr\"obner basis for an operad guarantees, as we saw above, that our operad is Koszul. Koszul duality for operads is used a lot for studying various operads, and proving Koszulness for particular operads is often an important and difficult problem, so an algorithm that provides a criterion of Koszulness can be very helpful.

\begin{example}
For the operad $\Lie$ of Lie algebras, the space of generators is one-dimensional. The element $\underline{[[a_1, a_2], a_3]}-[[a_1, a_3], a_2]-[a_1, [a_2,a_3]]$ forms a Gr\"obner basis of the ideal that defines $\Lie$. Indeed, there are no nontrivial reductions, the action of the symmetric group preserves the linear span of our element, and the $S$-polynomial 
 $$
S=\underline{[[[a_1,a_2],a_4],a_3]}+[[a_1,a_2],[a_3,a_4]]-[[[a_1,a_3],a_2],a_4]-[[a_1,[a_2,a_3]],a_4]
 $$
(it corresponds to the common multiple $[[[a_1,a_2],a_3],a_4]$) can be reduced to zero by the following sequence of reductions:
\begin{multline*}
S\mapsto 
[[[a_1,a_4],a_2],a_3]+[[a_1,a_2],[a_3,a_4]]-\underline{[[[a_1,a_3],a_2],a_4]}-\\
-[[a_1,[a_2,a_3]],a_4]+[[a_1,[a_2,a_4]],a_3]\mapsto 
-[[[a_1,a_3],[a_2,a_4]]+\underline{[[[a_1,a_4],a_2],a_3]}+\\ 
-[[[a_1,a_3],a_4],a_2]+[[a_1,a_2],[a_3,a_4]]-[[a_1,[a_2,a_3]],a_4]+[[a_1,[a_2,a_4]],a_3]
\mapsto\\
\mapsto-\underline{[[[a_1,a_3],a_4],a_2]}-[[a_1,a_3],[a_2,a_4]]+[[a_1,a_4],[a_2,a_3]]+[[[a_1,a_4],a_3],a_2]+\\
+[[a_1,a_2],[a_3,a_4]]-[[a_1,[a_2,a_3]],a_4]+[[a_1,[a_2,a_4]],a_3]
\mapsto-\underline{[[a_1,[a_2,a_3]],a_4]}-\\
-[[a_1,[a_3,a_4]],a_2]-[[a_1,a_3],[a_2,a_4]]+[[a_1,a_4],[a_2,a_3]]+[[a_1,a_2],[a_3,a_4]]+\\
+[[a_1,[a_2,a_4]],a_3]
\mapsto\underline{[[a_1,[a_2,a_4]],a_3]}-[a_1,[[a_2,a_3],a_4]]-[[a_1,[a_3,a_4]],a_2]-\\
-[[a_1,a_3],[a_2,a_4]]+[[a_1,a_2],[a_3,a_4]]\mapsto \underline{[[a_1,a_2],[a_3,a_4]]}+[a_1,[[a_2,a_4],a_3]]-\\
-[a_1,[[a_2,a_3],a_4]]-[[a_1,[a_3,a_4]],a_2]\mapsto[a_1,[a_2,[a_3,a_4]]]+[a_1,[[a_2,a_4],a_3]]-\\
-\underline{[a_1,[[a_2,a_3],a_4]]}\mapsto 0.
\end{multline*}
The basis of the operad $\Lie$ that we recover from this Gr\"obner basis, is, as it is easy to check, formed by all Lie monomials of the form $[A,B]$, where $A$ and $B$ are basis monomials of smaller degree, $A$ contains~$a_1$ and $B$ contains~$a_2$. One can check that the arity~$n$ part of this basis coincides with the multilinear part of the Shirshov basis~\cite{Ufn} of the free Lie algebra generated by $a_1,\ldots,a_n$ (for the ordering $a_1>\ldots>a_n$). It is also worth mentioning that this basis is essentially the same as the one that appeared recently in a work of Salvatore and Tauraso~\cite{ST}. 
\end{example}

\begin{example}
The operad $\As$ of associative algebras is the simplest example of an operad for which the space of generators is two-dimensional (it is one-dimensional if we consider this operad as a symmetric operad, but for the shuffle category we need two operations). Let us put 
 $$
\alpha(a_1,a_2)=a_1a_2,\quad \beta(a_1,a_2)=a_2a_1
 $$
(here $a,b\mapsto ab$ is the associative product which generates this operad as a symmetric operad).
Let us put $\alpha>\beta$. Then the ideal of relations in our operad is generated by the elements
\begin{gather*}
\underline{\alpha(\alpha(a_1,a_2),a_3)}-\alpha(a_1,\alpha(a_2,a_3)),\\
\underline{\alpha(\beta(a_1,a_2),a_3)}-\beta(\alpha(a_1,a_3),a_2),\\
\underline{\alpha(\alpha(a_1,a_3),a_2)}-\alpha(a_1,\beta(a_2,a_3)),\\
\underline{\alpha(\beta(a_1,a_3),a_2)}-\beta(\alpha(a_1,a_2),a_3),\\
\beta(a_1,\alpha(a_2,a_3))-\underline{\beta(\beta(a_1,a_3),a_2)},\\
\beta(a_1,\beta(a_2,a_3))-\underline{\beta(\beta(a_1,a_2),a_3)}.
\end{gather*}
These elements form a Gr\"obner basis: the action of the symmetric group preserves the subspace spanned by these elements, there are no reductions, and all the $S$-polynomials can be reduced to zero; for example, 
the $S$-polynomial
 $$
S=\underline{\beta(\alpha(\beta(a_1,a_3),a_4),a_2)}-\alpha(\beta(a_1,\alpha(a_2,a_3)),a_4)
 $$
(it corresponds to the common multiple $\alpha(\beta(\beta(a_1,a_3),a_2),a_4)$ of the leading terms of the second and the fifth basis elements), which can be reduced to zero by the following sequence of reductions:
\begin{multline*}
S\mapsto\underline{\beta(\beta(\alpha(a_1,a_4),a_3),a_2)}-\alpha(\beta(a_1,\alpha(a_2,a_3)),a_4)\mapsto\\ \mapsto
\beta(\alpha(a_1,a_4),\alpha(a_2,a_3))-\underline{\alpha(\beta(a_1,\alpha(a_2,a_3)),a_4)}\mapsto 0.
\end{multline*}
If we use the Gr\"obner basis to write down the basis for our operad, we shall see, for example, that the resulting basis for the space of ternary operations is 
 $$
a_1(a_2a_3), a_1(a_3a_2), a_2(a_1a_3), (a_2a_3)a_1, a_3(a_1a_2), (a_3a_2)a_1.
 $$
\end{example}

\begin{example}
The operads $\PreLie$ of pre-Lie algebras and its Koszul dual operad $\PreLie^!=\Perm$ were studied in several papers, see, for example, \cite{Chapoton,ChapotonPL}. As a symmetric operad, the operad $\PreLie$ is generated by one operation $\star$ which satisfies the relation
 $$
(a\star b)\star c-a\star(b\star c)=(a\star c)\star b-a\star(c\star b),
 $$
which guarantees that the bracket $[a,b]=a\star b-b\star a$ satisfies the Jacobi identity.

The operad $\Perm$ is generated by one operation~$\cdot$ which satisfies the relations
\begin{gather*}
(a\cdot b)\cdot c=a\cdot(b\cdot c),\\
a\cdot(b\cdot c)=a\cdot(c\cdot b).
\end{gather*}

These operads are PBW. For example, consider the operad $\PreLie$ as a shuffle operad with generators $\alpha\colon a,b\mapsto a\star b$ and $\beta\colon a,b\mapsto b\star a$. For the ordering $\alpha>\beta$, the quadratic relations are
\begin{gather*}
\underline{\alpha(\alpha(a_1,a_2),a_3)}-\alpha(a_1,\alpha(a_2,a_3))-
\alpha(\alpha(a_1,a_3),a_2)+\alpha(a_1,\beta(a_2,a_3)),\\
\underline{\alpha(\beta(a_1,a_2),a_3)}-\beta(\alpha(a_1,a_3),a_2)-
\beta(a_1,\alpha(a_2,a_3))+\beta(\beta(a_1,a_3),a_2),\\
\underline{\alpha(\beta(a_1,a_3),a_2)}-\beta(\alpha(a_1,a_2),a_3)-
\beta(a_1,\beta(a_2,a_3))+\beta(\beta(a_1,a_2),a_3),
\end{gather*}
and all $S$-polynomials for common multiples of the first leading term with the second and the third one (which themselves do not have nontrivial common multiples) can be reduced to zero. One can easily check that for the 
ordering $\alpha<\beta$ there exist $S$-polynomials that do not reduce to zero.
\end{example}

\begin{example}
The operad~$\Leib$ of Leibniz algebras~\cite{LodayLeib} is generated, as a symmetric operad, by the bracket $a_1,a_2\mapsto[a_1,a_2]$ (without any symmetries) such that the relation
 $$
[a,[b,c]]=[[a,b],c]-[[a,c],b]
 $$ 
is satisfied in every algebra over this operad. To interpret this operad as a shuffle operad, we introduce operations 
$\alpha(a,b)=[a,b]$ and $\beta(a,b)=[b,a]$. In terms of these operations, the defining relations for the shuffle operad $\Leib$ are
\begin{gather*}
\alpha(a_1,\alpha(a_2,a_3))-\alpha(\alpha(a_1,a_2),a_3)+\alpha(\alpha(a_1,a_3),a_2)=0,\\
\beta(a_1,\alpha(a_2,a_3))-\alpha(\beta(a_1,a_2),a_3)+\beta(\alpha(a_1,a_3),a_2)=0,\\
\beta(a_1,\beta(a_2,a_3))+\beta(\alpha(a_1,a_2),a_3)-\alpha(\beta(a_1,a_3),a_2)=0,\\
\alpha(a_1,\alpha(a_2,a_3))+\alpha(a_1,\beta(a_2,a_3))=0,\\
\beta(\alpha(a_1,a_2),a_3)+\beta(\beta(a_1,a_2),a_3)=0,\\
\beta(\alpha(a_1,a_3),a_2)+\beta(\beta(a_1,a_3),a_2)=0.
\end{gather*}

Let $\alpha>\beta$. Consider the version of the path-lexicographic ordering of tree monomials which compares words using the reverse degree-lexicographic ordering. 
It is straightforward to check that all $S$-polynomials can be reduced to zero, and the operad~$\Leib$ is PBW. 
\end{example}

\subsection{Operads from commutative algebras}

The following construction of an operad from a graded commutative algebra was introduced by the second author in~\cite{Khor}.

Let $A$ be a graded commutative algebra. Define an operad $\calO_A$ as follows. We put $\calO_A(n):=A_{n-1}$, and let
the partial composition map 
 $$
\circ_i\colon\calO_A(k)\otimes\calO_A(l)=A_{k-1}\otimes A_{l-1}\to A_{k+l-2}=\calO_A(k+l-1)  
 $$
be the product in~$A$. If the algebra $A$ is quadratic, then the operad $\calO_A$ is quadratic as well.

\begin{theorem}
If the algebra $A$ has a $k$-triangular basis (algebras are particular cases of operads, so our definition applies), then the operad $\calO_A$ has a $k$-triangular basis as well. 
\end{theorem}

\begin{proof}
Let us take some basis~$B$ of the algebra~$A$. It gives rise to a set of elements of~$\calO_A$ as follows. We define the collection $C$ of tree monomials in~$\calF_{V}$ (where~$V$ is the space of generators of~$A$) together with a one-to-one correspondence $\psi\colon C\to B$ recursively:
\begin{itemize}
 \item all generators belong to~$C$ and the correspondence~$\psi$ on them is tautological;
 \item if $\alpha\in C$ and $\beta\in C$, $\psi(\alpha)\psi(\beta)$ belongs to~$B$, and the arity of $\psi(\alpha)$ is equal to~$m$, then the element $\gamma:=\alpha\circ_m\beta$ belongs to~$C$, and $\psi(\gamma)=\psi(\alpha)\psi(\beta)$.
\end{itemize}
Then the image of $C$ under the projection from the free operad is a basis of $\calO_A$. Moreover, if $B$ was a $k$-triangular basis of $A$, then $C$ is a $k$-triangular basis of $\calO_A$. 
\end{proof}

In~\cite{Khor}, distributive lattices were used to prove that if the algebra $A$ is Koszul, then the operad $\calO_A$ is Koszul. Our results provide a somewhat simpler proof of that statement in the case when the algebra~$A$ is a PBW algebra.

\begin{definition}[\cite{DK}]
The operad $\Lie^2$ (called also the operad of two compatible brackets) is generated by two skew-symmetric operations (brackets) $\{\cdot,\cdot\}$ and $[\cdot,\cdot]$. The relations in this operad mean that all linear combinations of these brackets satisfy the Jacobi identity.
It is equivalent to the following identities in each algebra over this operad:
\begin{gather*}
\{a,\{b,c\}\}+\{b,\{c,a\}\}+\{c,\{a,b\}\}=0,\\
[a,\{b,c\}]+[b,\{c,a\}]+[c,\{a,b\}]+\{a,[b,c]\}+\{b[c,a]\}+\{c,[a,b]\}=0,\\
[a,[b,c]]+[b,[c,a]]+[c,[a,b]]=0.
\end{gather*}
The operad $\,{}^2\Com$ of two strongly compatible commutative products is generated by two symmetric binary operations (products) $\circ$ and $\bullet$ such that in any algebra over this operad the following identities hold:
\begin{gather*}
a\circ (b\circ c)=(a\circ b)\circ c,\\
a\circ(b\bullet c)=a\bullet(b\circ c)=b\circ(a\bullet c)=b\bullet(a\circ c)=c\circ(a\bullet b)=c\bullet(a\circ b),\\
a\bullet(b\bullet c)=(a\bullet b)\bullet c.
\end{gather*}
\end{definition}

The following result was proved in \cite{Strohmayer} (our original proof in \cite{DK} contained a gap). Our previous statement allows us to give yet another proof of this fact.

\begin{proposition}
The operad $\,{}^2\Com$ is Koszul. 
\end{proposition}

\begin{proof}
This operad is isomorphic to the operad $\calO_A$ for $A=\k[x,y]$.
\end{proof}

\begin{remark}
The operad $\Lie^k$ of $k$ compatible Lie brackets is generated by $k$ binary operations for which all linear combinations satisfy the Jacobi identity. It is easy to see that the Koszul dual operad is isomorphic to the operad $\calO_A$ for $A=\k[x_1,\ldots,x_k]$. Consequently, this operad is Koszul.
\end{remark}

\begin{definition}[\cite{HW}]
The operad $k-\Lie$ of Lie $k$-algebras is a quadratic operad with one $k$-ary skew-symmetric operation~$\omega$ satisfying the Jacobi identity 
 $$
\sum_{\sigma\in S_{2k-1}}(-1)^{\sigma}\omega(\omega(a_{\sigma(1)},\ldots,a_{\sigma(k)}),a_{\sigma(k+1)},\ldots,a_{\sigma(2k-1)})=0.
 $$
\end{definition}

\begin{proposition}
The operad $k-\Lie$ is Koszul. 
\end{proposition}

\begin{proof}
Its Koszul dual operad is, as it is easy to see, isomorphic to the operad obtained from the Koszul algebra $\k[t]$, where $\deg(t)=k-1$. This algebra is PBW, so the corresponding operad is PBW and hence Koszul. 
\end{proof}

\begin{definition}[\cite{Chap}]
The operad $\LieGriess$ is a binary quadratic operad with two skew-symmetric generators $[\cdot{,}\cdot]$ and $\{\cdot{,}\cdot\}$ that satisfy the identities
\begin{gather*}
\{a,\{b,c\}\}+\{b,\{c,a\}\}+\{c,\{a,b\}\}=0,\\
[a,\{b,c\}]+[b,\{c,a\}]+[c,\{a,b\}]+\{a,[b,c]\}+\{b,[c,a]\}+\{c,[a,b]\}=0,\\
\end{gather*}
\par
The Ramanujan operad $\Ram$ is a binary quadratic operad with a symmetric generator $\cdot\star\cdot$ and two skew-symmetric generator
$[\cdot{,}\cdot]$ and $\{\cdot{,}\cdot\}$ for which the product $\cdot\star\cdot$ generates a suboperad isomorphic to $\Com$, the operations $[\cdot{,}\cdot]$ and $\{\cdot{,}\cdot\}$ generate a suboperad isomorphic to $\LieGriess$, and these suboperads together are related by a distributive law~\cite{Markl,VD}
\begin{gather*}
[a,b\star c]=[a,b]\star c+b\star[a,c],\\
\{a,b\star c\}=\{a,b\}\star c+b\star\{a,c\}.
\end{gather*}
\end{definition}

\begin{corollary}
The operads $\LieGriess$ and $\Ram$ are Koszul.
\end{corollary}

\begin{proof}
One can easily check that the operad $\LieGriess^!$ is isomorphic to the operad $\calO_A$ with $A=\k[x,y]/(x^2)$. Thus this operad is Koszul, and
so is its dual $\LieGriess$. The distributive laws criterion implies that the operad~$\Ram$ is Koszul as well.
\end{proof}

This allows us to prove Chapoton's conjecture on the bigraded characters of the operad~$\Ram$. Recall that we can put the degree of $\cdot\star\cdot$ equal to $(0,0)$, the degree of $\{\cdot{,}\cdot\}$ equal to $(0,1)$, and the degree of $[\cdot{,}\cdot]$ equal to $(1,1)$, thus making the operad~$\Ram$ bigraded. Define the bigraded characters of its components, adding up dimensions of the homogeneous part of degree $(i,j)$ multiplied by $x^iy ^{j-i}$. The following result conjectured by Chapoton~\cite{Chap} can be deduced from the explicit description of the operad $\LieGriess$ above, together with the functional equation for dimensions and characters of an operad and its Koszul dual~\cite{GK}; details are left to the reader.
 
\begin{corollary}
The bigraded characters of the component $\Ram(n)$ of the operad~$\Ram$ is given by the
$n\text{th}$ Ramanujan polynomial~$\psi_n(x,y)$. Ramanujan polynomials can be defined recursively
as follows:
\begin{align*}
  \psi_1&=1,\\
  \psi_{n+1}&=\psi_n + (x+y) (n \psi_n +x \partial_x \psi_n). 
\end{align*}
\end{corollary}

\section{Further directions}\label{further} 

Our results can be generalised to the case of dioperads~\cite{Gan}, $\frac12$PROPs~\cite{MV}, and coloured operads~\cite{vdL} in a rather straightforward way. For example, the dioperad of Lie bialgebras has a quadratic Gr\"obner 
basis; it gives yet another proof of its Koszulness. We shall discuss the details and more examples elsewhere. Currently, we do not know whether Gr\"obner bases machinery for properads and PROPs can be defined in a similar 
fashion.

Shuffle operads can be used for questions of homological algebra as well: for a symmetric operad, its symmetric (co)bar construction coincides with its shuffle (co)bar construction, so we can compute its cohomology using the 
shuffle category. It turns out that in that category one can use a Gr\"obner basis to construct a remarkable free resolution of the trivial module. We shall describe it in our next paper.

Our analogue of Buchberger's algorithm should be interpreted as a call for a computer algebra system that would compute Gr\"obner bases for operads presented by generators and relations, and use those bases to compute dimensions 
of graded components, normal forms of elements etc. A first step in this direction is a \texttt{Haskell} package \texttt{Operads}~\cite{VD+MVJ}.

For some particular operads which has not been studied thoroughly yet, Gr\"obner bases might be the right tool to approach them with. In particular, they might be quite useful in the case of operads for which the generating 
operations are not binary. An interesting examples of such operads are the operad of totally and partially associative $k$-ary algebras from~\cite{Gned} (see an interesting recent paper~\cite{MR} for a detailed study of this operad). 
An algebraic structure with a ternary generating operation for which the corresponding operad has been neglected so far controls generalised 3-Lie algebras of Cherkis--Saemann~\cite{CS}.

\end{document}